\documentclass[10pt,a4paper]{article}

\usepackage{amsmath}
\usepackage{amsthm}
\usepackage{amssymb}
\usepackage{amsfonts}

\newtheorem*{icet}{$I$CET}
\newtheorem*{gcp}{GCP}
\newtheorem*{mucet}{$\mu$CET}
\newtheorem*{cwlln}{CWLLN}
\newtheorem*{mucwlln}{$\mu$CWLLN}
\newtheorem*{egcp}{EGCP}
\newtheorem*{bcp}{BCP}

\newtheorem*{mm}{MaxProb/MaxEnt}
\newtheorem*{mugcp}{$\mu$GCP}
\newtheorem*{ST}{Sanov's Theorem}
\newtheorem*{lem}{Lemma}
\newtheorem*{rCWLLN}{($r=2$)-tuple CWLLN}
\newtheorem{lemN}{Lemma}

\theoremstyle{definition}
\newtheorem{ex}{Example}

\theoremstyle{remark}

\begin{document}
\title{%
Conditional Equi-concentration of Types\thanks{This work was
supported by VEGA 1/0264/03 grant. Valuable discussions with George
Judge, Brian R. La Cour, Alberto Solana-Ortega, Ondrej \v Such and
Viktor Witkovsk\' y are gratefully acknowledged. Lapses are mine.}}
\author{
M. Grendar\thanks{E-mail address: marian.grendar@savba.sk}\\
Department of Mathematics, FPV UMB, \\ Tajovskeho 40, 974 01 Banska
Bystrica, Slovakia, \\
Institute of Mathematics and CS, Slovak Academy of Sciences (SAS) \\
and Institute of Measurement Science of SAS\\
\\
\small{To Mar, in memoriam}}
\date{}
\maketitle

\begin{abstract}
Conditional Equi-concentration of Types on $I$-projections  is
presented. It provides an extension of Conditional Weak Law of Large
Numbers to the case of several $I$-projections. Also a multiple
$I$-projections  extension of Gibbs Conditioning Principle is
developed. $\mu$-projection variants of the probabilistic laws are
stated. Implications of the results for Relative Entropy
Maximization, Maximum Probability, Maximum Entropy in the Mean and
Maximum R\'enyi-Tsallis Entropy methods are discussed.
\end{abstract}


\vspace{.08in}\noindent \textbf{2000 Mathematics Subject
Classification} Primary 62B15; Secondary 94A15, 62J20

\section{Terminology and Notation}

Let $\{X\}_{l=1}^n$ be a sequence of independently and identically
distributed random variables with a common law (measure)  on a
measurable space. Let the measure be concentrated on finite number
$m$ of atoms from the set $\mathrm{X} = \{x_1, x_2, \dots, x_m\}$
called support or alphabet.  Let $q_i$ denote the probability
(measure) of $i$-th element of $\mathrm X$; $q$ will be assumed
strictly positive and called source or generator. Let
$\mathrm{P}(\mathrm{X})$ be a set of all probability mass functions
(pmf's) on $\mathrm X$.

A type  (also called $n$-type, empirical measure, frequency
distribution or occurrence vector) induced by a sequence
$\{X\}_{l=1}^n$ is  pmf $\nu^n \in \mathrm{P}(\mathrm{X})$ whose
$i$-th element $\nu_i^n$ is defined as: $\nu_i^n = n_i/n$ where $n_i
= \sum_{l=1}^n I(X_l = x_i)$; there $I(\cdot)$ is the indicator
function. Multiplicity $\Gamma(\nu^n)$ of type $\nu^n$ is:
$\Gamma(\nu^n) = n!/ \prod_{i=1}^m n_i!$.

Let $\mathrm{\Pi} \subseteq \mathrm{P}(\mathrm{X})$. Let
$\mathrm{P}_n$ denote a subset of $\mathrm{P}(\mathrm{X})$ which
consists of all $n$-types. Let $\mathrm{\Pi}_n = \mathrm\Pi \cap
\mathrm{P}_n$.

On $\mathrm P(\mathrm X)$ topology induced by the standard topology
on $\mathrm{R}^m$ is assumed.

$\mu$-projection $\hat{\nu}^n$ of $q$ on $\mathrm{\Pi}_n \neq
\emptyset$ is defined as: $\hat{\nu}^n = \arg \, \sup_{\nu^n \in
\mathrm{\Pi}_n} \pi(\nu^n; q)$, where $\pi(\nu^n; q) = \Gamma(\nu^n)
\prod (q_i)^{n \nu^n_i}$. Alternatively, the $\mu$-projection can be
defined as $\hat{\nu}^n = \arg \, \sup_{\nu^n \in \mathrm{\Pi}_n}
\pi(\nu^n | \nu^n \in \mathrm{\Pi}_n; q)$, where $\pi(\nu^n| \nu^n
\in \mathrm{\Pi}_n; q)$ denotes the conditional probability that if
an $n$-type belongs to $\mathrm{\Pi}_n$ then it is just the type
$\nu^n$. $\mu$-projection can be also equivalently defined as the
supremum of  posterior probability,  cf. \cite{gg_bayes}.

$I$-projection $\hat p$ of $q$ on $\mathrm\Pi$ is  $\hat p = \arg \,
\inf_{p \in \mathrm\Pi} I(p || q)$, where $I(p || q) =
\sum_{\mathrm{X}} p_i \log\frac{p_i}{q_i}$ is the $I$-divergence
(also known as Kullback-Leibler distance, $\pm$ relative entropy).

$\pi(\nu^n \in \mathrm{A}| \nu^n \in \mathrm{B}; q)$ will denote the
conditional probability that if a type drawn from $q \in
\mathrm{P}({\mathrm X})$ belongs to $\mathrm{B} \subseteq
\mathrm{\Pi}$ then it belongs to $\mathrm{A} \subseteq
\mathrm{\Pi}$.

\section{Boltzmann Jaynes Inverse Problem and Conditional Law of Large Numbers}

Having the terminology introduced, Boltzmann Jaynes Inverse Problem
(BJIP) can be stated as follows: there is the source $q$ and a set
$\Pi_n$ of $n$-types. It is necessary to select an $n$-type (one or
more) from the set $\Pi_n$. To solve BJIP it is necessary to provide
an algorithm for selection of type from $\Pi_n$ when the
information-quadruple $\{\mathrm{X}, n, q, \mathrm{\Pi}_n\}$ and
nothing else is supplied. Clearly, if $\Pi_n$ contains more than one
type, BJIP becomes an under-determined and in this sense ill-posed
problem.

Usually, BJIP is solved by means of the method of Relative Entropy
Maximization (REM/MaxEnt). This is mostly done  for  $n \rightarrow
\infty$. In this case the set of types $\Pi_n$ effectively turns
into a set of probability mass functions $\Pi$.

Typically,  $\mathrm\Pi$ is defined by moment consistency
constraints (mcc) of the following form\footnote{In the simplest
case of single non-trivial constraint.}: $\mathrm\Pi_{mcc} = \{p:
\sum_{i=1}^m p_i u_i = a, \sum_{i=1}^m p_i = 1\}$, where $a \in
\mathrm{R}$ is a given number, $u$ is a given vector. The feasible
set $\mathrm\Pi_{mcc}$ which mcc define is convex and closed.
$I$-projection $\hat{p}$ of $q$ on $\mathrm\Pi_{mcc}$ is unique and
belongs to the exponential family of distributions; $\hat{p}_i =
k(\lambda) q_i e^{-\lambda u_i}$, where $k(\lambda) = 1/\sum_{i=1}^m
q_i e^{-\lambda u_i}$, and $\lambda$ is such that $\hat{p}$
satisfies mcc.

In the case of BJIP with  $\mathrm{\Pi}_{mcc}$, or in general for
any closed, convex, rare set $\mathrm\Pi$, application of REM/MaxEnt
method is justified by Conditional Weak Law of Large Numbers
(CWLLN). CWLLN, in its textbook form, reads \cite{CT}:

\begin{cwlln}
 Let $\mathrm{X}$ be a finite set. Let $\mathrm{\Pi}$ be
a closed, convex set which does not contain $q$. Let $n \rightarrow
\infty$. Then for $\epsilon
> 0$ and $i = 1, 2, \dots, m$,
$$
\lim_{n \rightarrow \infty} \pi(|\nu^n_i - \hat{p}_i| \le \epsilon
\,|\, \nu^n \in \mathrm{\Pi}; q) = 1.
$$
\end{cwlln}

CWLLN says that if types  are confined to the set $\mathrm\Pi$ then
they asymptotically conditionally concentrate on the $I$-projection
$\hat{p}$ of the source of types $q$ on the set $\mathrm\Pi$.
Stated, informally, from another side: if a source $q$ is confined
to produce types from convex and closed $\mathrm{\Pi}$ it is
asymptotically conditionally 'almost impossible' to find a type
other than the one which has the highest/supremal value of relative
entropy with respect to $q$.

Conditional Weak Law of Large Numbers emerged from a series of works
which include  \cite{Sanov}, \cite{Bartfai}, \cite{Vincze},
\cite{L},  \cite{GOR}, \cite{Vasicek}, \cite{Zabel}, \cite{CC},
\cite{Topsoe},  \cite{Csi}, \cite{Ellis}, \cite{DZ}, \cite{LPS},
\cite{LN}. For new developments see \cite{H}.

An information-theoretic  proof (see \cite{CT}) of CWLLN utilizes
so-called Py\-tha\-go\-re\-an theorem (cf. \cite{Chentsov}),
Pin\-sker inequality and standard inequalities for factorial. The
Py\-tha\-go\-re\-an theorem is known to hold for closed convex sets.
Alternatively, CWLLN can be obtained as a consequence of Sanov's
Theorem (ST). The ST-based proof of CWLLN will be recalled here.
First, Sanov's Theorem and its proof (adapted from \cite{CsiMT},
\cite{CT}).

\begin{ST}
 Let $\mathrm X$ be finite. Let $\mathrm A \subseteq \mathrm\Pi$
be an open set. Then
 \begin{equation*}
  \lim_{n \rightarrow \infty} \frac{1}{n} \log \pi(\nu^n \in \mathrm A) = -
  I(\hat{p} || q),
 \end{equation*}
where $\hat{p}$ is an $I$-projection of $q$ on $\mathrm A$.
\end{ST}

\begin{proof} \ \ \cite{CT}, \cite{CsiMT} \ \
$\pi(\nu^n \in \mathrm A)  = \sum_{\nu^n \in \mathrm A} \pi(\nu^n;
q)$. Upper and lower bounds on $\pi(\nu^n; q)$ (recall proof of the
Lemma at Appendix):
\begin{equation*}
\left(\frac{m}{n}\right)^m \prod_{i=1}^m
\left(\frac{q_i}{{\nu}_i^n}\right)^{n {\nu}_i^n} < \pi(\nu^n; q) \le
\prod_{i=1}^m \left(\frac{q_i}{\nu_i^n}\right)^{n \nu_i^n}.
\end{equation*}

$\sum_{\nu^n \in \mathrm A} \pi(\nu^n; q) < N \prod_{i=1}^m
(\frac{q_i}{\hat{\nu}_i^n})^{n \hat{\nu}_i^n}$, where $N$ stands for
number of all $n$-types and $\hat{\nu}_i^n$ is an $I$-projection of
$q$ on $\mathrm{A}_n = \mathrm A \cap \mathrm{P}_n$ (i.e., any of
the $n$-types which attain supremal value of $\prod_{i=1}^m
(\frac{q_i}{\hat{\nu}_i^n})^{n \hat{\nu}_i^n}$). $N$ is smaller than
$(n + 1)^m$.

Thus
\begin{multline*}
\frac{1}{n}  \left(n\sum_{i=1}^m \hat{\nu}_i^n
\log\frac{q_i}{\hat{\nu}_i^n} + m(\log m - \log n)\right) <
\frac{1}{n} \log \pi(\nu^n \in \mathrm A) \\ < \frac{1}{n}
\left(n\sum_{i=1}^m \hat{\nu}_i^n \log\frac{q_i}{\hat{\nu}_i^n} +  m
\log(n+1)\right).
\end{multline*}

Since $\mathrm A$ is by the assumption open and under the maintained
assumption of strictly positive $q$ it is also continuous, $\lim_{n
\rightarrow \infty} \sum \hat{\nu}_i^n \log\frac{q_i}{\hat{\nu}_i^n}
= \sum \hat{p}_i \log\frac{q_i}{\hat{p}_i}$, where $\hat{p}$ is an
$I$-projection of $q$ on $\mathrm A$. Thus, for $n \rightarrow
\infty$ the upper and lower bounds on $\frac{1}{n} \log \pi(\nu^n
\in \mathrm A)$ collapse into $\sum_{i=1}^m \hat{p}_i
\log\frac{q_i}{\hat{p}_i}$.
\end{proof}

\begin{proof}[A proof of CWLLN]\ \ \cite{CsiMT} \ \
Let $\mathrm A = \{p: |p_i - \hat{p}_i| > \epsilon, i = 1, 2, \dots,
m\}$. Then ST can be applied to it, leading $\lim_{n \rightarrow
\infty} \frac{1}{n} \log \pi(\nu^n \in \mathrm A | \nu^n \in
\mathrm\Pi; q) = - (I(\hat{p}_{\mathrm A} || q) -
I(\hat{p}_{\mathrm\Pi} || q))$. Since $I(\hat{p}_{\mathrm A} || q) -
I(\hat{p}_{\mathrm\Pi} || q) > 0$ and since the set $\mathrm\Pi$
admits unique $I$-projection (the uniqueness arises from the fact
that the set is convex and closed, and $I(\cdot || \cdot)$ is
convex), the proof is complete.
\end{proof}

CWLLN can be viewed as a special case of a stronger result, which is
commonly known as Gibbs Conditioning Principle (GCP), see Sect. 5.

\section{Motivation and Programme}

Frequency moment constraints considered by physicists (see for
instance \cite{RAD}) define a non-convex feasible set of probability
distributions which in general can admit multiple $I$-projections.
This work builds upon \cite{gg_wawa}, \cite{LS}, \cite{gg_aei},
\cite{gg_nonlin}, \cite{gg_asy} and aims to develop an extension of
CWLLN and Gibbs Conditioning Principle  to the case of multiple
$I$-projections.

It has also another goal: to introduce concept of $\mu$-projection
and to formulate $\mu$-projection variants of the probabilistic
laws. They, among other things, allow for a more elementary reading
than their $I$-projection counterparts. At the same time they
provide a probabilistic justification of Maximum Probability method
\cite{gg_what}.

The paper is organized as follows: in the next section some basic
questions regarding asymptotic behavior of conditional probability
are posed. Two illustrative examples are then used to introduce
Conditional Equi-con\-cen\-trati\-on of Types on $I$-projections.
Next, an extension of Gibbs Conditioning Principle - the stronger
form of CWLLN - is provided. Asymptotic identity of $I$-projections
and $\mu$-projections is discussed in Section 6 and $\mu$-variants
of the probabilistic laws are presented afterwards. Implications of
the results for Maximum Entropy, Maximum Probability and Maximum
R\'enyi-Tsallis Entropy methods are drawn at Section 8. Section 9
mentions in passing other related results: $r$-tuple extension of
CWLLN and Bayesian Conditional Law of Large Numbers. Section 10
summarizes the paper. Appendix contains a sketch of proof of $I$CET
and of Extended GCP. It also shows that concentration of types can
in some sense happen also on isolated $I$-projections, provided that
they are rational.

\section{Conditional Equi-concentration of Types}

What happens when $\mathrm\Pi$ admits multiple $I$-projections? Do
the conditional concentration of types happens on them? If yes, do
the type concentrate on each of them? If yes, what is the
proportion? In order to address these questions, it is instrumental
to consider a couple of examples.

\begin{ex} \cite{gg_aei}  \ \  Let $\mathrm{\Pi} = \{ p: \sum_{i=1}^m
p_i^\alpha - a = 0, \sum_{i=1}^m p_i - 1 = 0\}$, where $\alpha, a
\in \mathrm{R}$. Note that the first constraint, known as frequency
constraint, is non-linear in $ p$ and $\mathrm{\Pi}$ is for
$|\alpha| > 1$ non-convex.

Let $\alpha = 2$, $m = 3$ and $a = 0.42$ (the value was obtained for
$p = [0.5\ 0.4\ 0.1]$). Then there are the following three
$I$-projections of uniform distribution $ q = [1/3 \ 1/3 \ 1/3]$ on
$\mathrm{\Pi}$: $\hat{ p}_1 = [0.5737 \ 0.2131 \ 0.2131]$, $\hat{
p}_2 = [0.2131 \ 0.5737$ \  $\ 0.2131]$ and $\hat{ p}_3 = [0.2131 \
0.2131$ $\ 0.5737]$
 (see \cite{gg_nonlin}). Note that they form a group of permutations. As it will become clear
 later, it suffices to investigate convergence to say $\hat{ p}_1$.

For $n=30$ there are only two groups of types in $\mathrm{\Pi}$: G1
comprises $[0.5666\ 0.2666$ $\ 0.1666]$ and five other permutations;
G2 consists of $[0.5\ 0.4\ 0.1]$ and the other five permutations.
So, together there are $12$ types.

Value of the square of the Euclidean distance $\delta$ between $\nu$
and $\hat{ p}_1$ attains its minimum $\delta_{G1} = 0.0051$ within
G1 group for the following two types: $[0.5666$ $0.2666 \  \
0.1666]$, $[0.5666\ 0.1666\ 0.2666]$. Within G2 the smallest
$\delta_{G2} = 0.0532$ is attained by $[0.5\ 0.4\ 0.1]$ and $[0.5\
0.1\ 0.4]$.

Thus, if $\epsilon = \epsilon_1$ is chosen so that the
$\epsilon$-ball $B(\hat{ p}_1, \epsilon_1)$ centered at $\hat{p}_1$
contains only the two types from G1 (which at the same time
guarantees that $\hat{ p}_1$ is the only $I$-projection in the
ball), then $\mathbb{\pi}(\nu \in B(\hat{ p}_1, \epsilon_1) | \nu
\in \mathrm{\Pi}) = 2 \cdot 0.1152 = 0.2304$. In words: probability
that if $ q$ generated a type from $\mathrm{\Pi}$ than the type
falls into the ball containing only types which are closest to the
$I$-projection is 0.2304. If $\epsilon = \epsilon_2$ is chosen so
that also the two types from G2 are included in the ball and also so
that it is the only $I$-projection in the ball (any $\epsilon_2 \in
(\sqrt{0.0532}, \sqrt{0.1253})$ satisfies both the requirements),
then $\mathbb{\pi}(\nu^n \in B(\hat{ p}_1, \epsilon_2) | \nu^n \in
\mathrm{\Pi}) = \frac{1}{3}$.

For $n=330$ there are four groups of types in $\mathrm{\Pi}$: G1, G2
and a couple of new one: G3 consists of $[0.4727\  0.4333\  0.0939]$
and all its permutations; G4  comprises the type $[0.5727\  0.2333\
0.1939]$ and its permutations. Hence, the total number of types from
$\mathrm{\Pi}$
 which are supported by random sequences of size
$n=330$ is  $24$.

$\delta_{G3}$ for the two types from G3 which are closest to $\hat{
p}_1$ is $0.0729$. The smallest $\delta_{G4} = 0.00077$ is attained
by $[0.5727\ 0.2333\  0.1939]$ and by $[0.5727\ 0.1939\  0.2333]$.
Thus, clearly, the two types from G4 have the smallest Euclidean
distance to $\hat{ p}_1$ among all types from $\mathrm{\Pi}$ which
are based on samples of size $n=330$. Again, setting $\epsilon$ such
that the ball $B(\hat{ p}_1, \epsilon)$ contains only the two  types
which are closest to $\hat{ p}_1$ leads to the $0.261$ value of the
conditional probability. Note the important fact, that the
probability has risen, as compared to the corresponding value 0.2304
for $n=30$.

Moreover,  if $\epsilon$ is set such that besides  the two types
from G4 also the second closest types (i.e. the two types from G1)
are included in the ball then  the conditional probability is
indistinguishable from $\frac{1}{3}$. Hence, there is virtually no
conditional chance of observing any of the remaining 4 types. The
same holds for the types which concentrate around $\hat{ p}_2$ or
$\hat{ p}_3$. Thus,  in total, a half of the 24 types is almost
impossible to observe.

The Example illustrates, that  the conditional probability of
finding a type which is close (in the Euclidean distance) to one of
the three $I$-projections goes to $\frac{1}{3}$.
\end{ex}

\smallskip

\begin{ex} \cite{gg_aei} \ \
 Let $\mathrm{\Pi} = \mathrm{\Pi}_1 \cup \mathrm{\Pi}_2$,
where $\mathrm{\Pi}_j = \{ p: \sum_{i=1}^m p_i x_i = a_j;$ $
\sum_{i=1}^m p_i = 1\}, j = 1, 2$. Thus $\mathrm{\Pi}$ is union of
two sets, each of whose  is given by the  moment consistency
constraint. If $ q$ is chosen to be the uniform distribution, then
values $a_1$, $a_2$ such that there will be two different
$I$-projections of the uniform $ q$ on $\mathrm{\Pi}$ with the same
value of $I$-divergence (as well as of the Shannon's entropy) can be
easily found. Indeed, for any $a_1 = \mu + \Delta$, $a_2 = \mu -
\Delta$, where $\mu = EX$ and $\Delta \in (0, (X_{\max} -
X_{\min})/2)$,  $\hat{ p}_1$ is just a permutation of $\hat{ p}_2$,
and as such attains the same value of Shannon's entropy. To see that
types which are based on random samples of size $n$ from
$\mathrm{\Pi}$ indeed concentrate on the $I$-projections with equal
measure note, that for any $n$ to each type in $\mathrm{\Pi}_1$
corresponds  a unique permutation of the type in $\mathrm{\Pi}_2$.
Thus, types in $\epsilon$-ball with center at $\hat{ p}_1$ have the
same conditional probabilities $\mathbb{\pi}$ as types in the
$\epsilon$-ball centered at $\hat{ p}_2$. This, together with
convexity and closed-ness of both $\mathrm{\Pi}_j$, for which the
conditional concentration of types on the respective $I$-projection
is established by CWLLN, directly implies that
$$
 \lim_{n \rightarrow \infty} \mathbb{\pi}(\nu \in B(\hat{ p}_j, \epsilon)
  | \nu^n \in \mathrm{\Pi}) = \frac{1}{2}    \qquad\text{$j = 1,
  2$}.
$$
\end{ex}

\smallskip

Conditional Equi-concentration of Types on $I$-projections ($I$CET)
attempts to capture behavior of the conditional measure which the
above Examples illustrate. To this end, notion of the proper
$I$-projection will be needed.

$I$-projection $\hat{p}$ of $q$ on $\mathrm\Pi$ will be called
proper if $\hat{p}$ is not  isolated point of $\mathrm\Pi$.

\smallskip

\begin{icet}
 Let $\mathrm{X}$ be finite. Let there be $\mathrm k$ proper
$I$-projections $\hat{p}^1, \hat{p}^2, \dots, \hat{p}^\mathrm{k}$ of
$q$ on $\mathrm\Pi$. Let $\epsilon > 0$ be such that for $j = 1, 2,
\dots, \mathrm{k}$ $\hat{p}^j$ is the only proper $I$-projection of
$q$ on $\mathrm\Pi$ in the ball $B(\hat{p}^j, \epsilon)$. Let $n
\rightarrow \infty$. Then for $j = 1, 2, \dots, \mathrm{k}$,
\begin{equation*}
\pi(\nu^n  \in B(\epsilon, \hat{p}^j) | \nu^n \in \mathrm{\Pi}; q) =
1/\mathrm{k}.
\end{equation*}
\end{icet}

$I$CET says, informally, that source/generator $q$, when confined to
produce types from a set $\mathrm{\Pi}$, - as $n$ gets large - hides
itself behind any of the proper $I$-projections equally likely.

Expressed in Statistical Physics terminology $I$CET says that each
of equilibrium points ($I$-projections) is asymptotically
conditionally equally possible. The Conditional Equi-concentration
of Types 'phenomenon' resembles the triple point phenomenon of
Thermodynamics.

A sketch of proof of $I$CET is relegated to the Appendix.

\section{Gibbs Conditioning Principle and its Extension}

Gibbs conditioning principle (cf. \cite{Csi}, \cite{DZ}, \cite{LN})
- also known as the stronger form of CWLLN - complements CWLLN by
stating that:

\begin{gcp}
 Let $\mathrm{X}$ be a finite set. Let $\mathrm{\Pi}$ be closed,
convex set. Let $n \rightarrow \infty$. Then for a fixed $t$,
$$
\lim_{n \rightarrow \infty} \pi(X_1 =x_1, \dots, X_t = x_t | \nu^n
\in \mathrm{\Pi}; q) = \prod_{l=1}^t \hat{p}_{x_l}.
$$
\end{gcp}

GCP, says, very informally, that if the source $q$ is confined to
produce sequences which lead to types in a set $\mathrm\Pi$ then
elements of any such  sequence (of fixed length $t$) behave
asymptotically conditionally as if they were drawn identically and
independently from the $I$-projection of $q$ on $\mathrm\Pi$ -
provided that the last is unique (among other things).

GCP was developed at \cite{Csi} under the name of conditional
quasi-independence of outcomes. Later on, it was brought into more
abstract form in large deviations literature, where it also obtained
the GCP name (cf. \cite{DZ}, \cite{LN}). A simple proof of GCP can
be found at \cite{CsiMT}. GCP is proven also for continuous alphabet
(cf. \cite{GOR},  \cite{CsiMT}, \cite{DZ}).

The following theorem provides an extension of GCP to the case of
multiple $I$-pro\-jec\-ti\-ons.

\begin{egcp}
 Let there be $\mathrm k$ proper $I$-projections $\hat{p}^1,
\hat{p}^2, \dots, \hat{p}^\mathrm{k}$ of $q$ on $\mathrm\Pi$. Then
for a fixed $t$ and $n \rightarrow \infty$,
\begin{equation*}
\pi(X_1 = x_1, \dots, X_t = x_t | \nu^n \in \mathrm\Pi; q) =
\frac{1}{\mathrm{k}} \sum_{j=1}^{\mathrm{k}} \prod_{l=1}^t
\hat{p}_{x_l}^j.
\end{equation*}
\end{egcp}

For $t=1$ Extended Gibbs Conditioning Principle (EGCP) says that the
conditional probability of a letter is asymptotically given by the
equal-weight mixture of proper $I$-projection probabilities of the
letter. For a general sequence, EGCP states that the conditional
probability of a sequence is asymptotically equal to the mixture of
joint probability distributions. Any  ($j$-th) of the $\mathrm{k}$
joint distributions is such as if the sequence was iid distributed
according to a ($j$-th)  proper $I$-projection.

A proof of EGCP is sketched at the Appendix.

\section{Asymptotic~Identity~of~$\mu$-Projections and $I$-Projections}

At (\cite{gg_what}, Thm 1 and its Corollary, aka MaxProb/MaxEnt Thm)
it was shown that maximum probability type converges to
$I$-projection; provided that $\mathrm{\Pi}$ is defined by a
differentiable constraints. A more general result which states
asymptotic identity of $\mu$-projections and $I$-projections for
general set $\mathrm{\Pi}$ was presented at \cite{gg_asy}.

\begin{mm}
 Let $\mathrm{X}$ be finite set. Let $\mathrm{M}_n$ be
set of all $\mu$-projections of $q$ on $\mathrm{\Pi}_n$. Let
$\mathrm{I}$ be  set of all $I$-projections of $q$ on
$\mathrm{\Pi}$. For $n \rightarrow \infty$, $\mathrm{M}_n =
\mathrm{I}$.
\end{mm}

Since $\pi(\nu^n; q)$ is defined for $\nu^n \in \mathrm{Q}^m$,
$\mu$-projection can be defined only for $\mathrm{\Pi}_n$ when $n$
is finite. The Thm permits to define a $\mu$-projection $\hat{\nu}$
also on $\mathrm{\Pi}$: $\hat{\nu} = \arg \sup_{r \in \mathrm{\Pi}}
- \sum_{i=1}^m r_i \log\frac{r_i}{q_i}$. The $\mu$-projections of
$q$ on $\mathrm\Pi$ and $I$-projections of $q$ on the same set
$\mathrm\Pi$ are undistinguishable.

It is worth highlighting that for a finite $n$, $\mu$-projections
and $I$-projections of $q$ on $\mathrm{\Pi}_n$ are in general
different. This explains why $\mu$-form of the probabilistic laws
deserves to be stated separately of the $I$-form; though formally
they are undistinguishable. Thus, the MaxProb/MaxEnt Thm (in its new
and to a smaller extent also in its old version) permits directly to
state $\mu$-projection variants of CWLLN, GCP, $I$CET and EGCP:
$\mu$CWLLN, $\mu$GCP, $\mu$CET and Boltzmann Conditioning Principle
(BCP).

\smallskip

\section{$\mu$-Variants of the Probabilistic Laws}

$\mu$-variant of CWLLN reads:

\begin{mucwlln}
 Let $\mathrm{X}$ be a finite set. Let $\mathrm{\Pi}$ be closed,
convex set. Let $n \rightarrow \infty$. Then for $\epsilon
> 0$ and $i = 1, 2, \dots, m$,
$$
\lim_{n \rightarrow \infty} \pi(|\nu^n_i - \hat{\nu}_i| < \epsilon |
\nu^n \in \mathrm{\Pi}; q) = 1.
$$
\end{mucwlln}

Core of $\mu$CWLLN can be loosely expressed as: {\it types, when
confined to a set} $\mathrm{\Pi}$, {\it conditionally concentrate on
the asymptotically most probable type} $\hat{\nu}$.

$\mu$-projection $\hat{\nu}$ of $q$ on $\mathrm\Pi$ will be called
proper if $\hat{\nu}$ is not  isolated point of $\mathrm\Pi$.

\begin{mucet}
 Let $\mathrm{X}$ be finite. Let there be $\mathrm k$
proper $\mu$-projections $\hat{\nu}^1, \hat{\nu}^2, \dots,
\hat{\nu}^\mathrm{k}$ of $q$ on $\mathrm\Pi$. Let $\epsilon > 0$ be
such that for $j = 1, 2, \dots, \mathrm{k}$ $\hat{\nu}^j$ is the
only proper $\mu$-projection of $q$ on $\mathrm\Pi$ in the ball
$B(\hat{\nu}^j, \epsilon)$. Let $n \rightarrow \infty$. Then for $j
= 1, 2, \dots, \mathrm{k}$,
\begin{equation*}
\pi(\nu^n  \in B(\epsilon, \hat{\nu}^j) | \nu^n \in \mathrm{\Pi}; q)
= 1/\mathrm{k}.
\end{equation*}
\end{mucet}

Core of $\mu$-variant of the Conditional Equi-concentration of Types
states, loosely, that types conditionally concentrate on each of the
asymptotically most probable types in equal measure.

\begin{mugcp}
 Let $\mathrm{X}$ be a finite set. Let $\mathrm{\Pi}$ be
closed, convex set. Let $n \rightarrow \infty$. Then for a fixed
$t$,
$$
\lim_{n \rightarrow \infty} \pi(X_1 =x_1, \dots, X_t = x_t | \nu^n
\in \mathrm{\Pi}; q) = \prod_{l=1}^t \hat{\nu}_{x_l}.
$$
\end{mugcp}

$\mu$-variant of EGCP deserves a special name. It will be called
Boltzmann Conditioning Principle (BCP).

\begin{bcp}
 Let there be $\mathrm k$ proper $\mu$-projections
$\hat{\nu}^1, \hat{\nu}^2, \dots, \hat{\nu}^\mathrm{k}$ of $q$ on
$\mathrm\Pi$. Then for a fixed $t$ and $n \rightarrow \infty$,
\begin{equation*}
\pi(X_1 = x_1, \dots, X_t = x_t | \nu^n \in \mathrm\Pi; q) =
\frac{1}{\mathrm{k}} \sum_{j=1}^{\mathrm{k}} \prod_{l=1}^t
\hat{\nu}_{x_l}^j.
\end{equation*}
\end{bcp}

\section{Implications}

The results have some  implications for application of REM, MaxProb
and Maximum R\'enyi-Tsallis Entropy methods to Boltzmann Jaynes
Inverse Problem.

\subsection{$I$- or $\mu$-Projection? MaxEnt or MaxProb?}

With $\mu$-projection Maximum Probability method (MaxProb,
\cite{gg_what}) is associated. Given the BJIP
in\-for\-mati\-on-quadruple $\{\mathrm{X}, n, q, \mathrm{\Pi}_n \}$,
MaxProb prescribes to select from $\mathrm{\Pi}_n$ type(s) which has
the supremal/maximal probability $\pi(\nu^n; q)$\footnote{A
technique for determination of $\mu$-projections was suggested at
\cite{gg_det}.}.

$\mu$-projections and $I$-pro\-jec\-tion\-s are asymptotically
indistinguishable. In plain words: for $n \rightarrow \infty$ the
Relative Entropy Maximization method (REM/Max\-Ent) (either in its
Jaynes' \cite{Jaynes}, \cite{Jbook} or Csisz\'ar's interpretation
\cite{CsiME}) selects the same distribution(s) as MaxProb (in its
more general form which instead of the maximum probable types
selects supremum-probable $\mu$-projections). This result (in the
older form, \cite{gg_what}) was at \cite{gg_what} {\it interpreted}
as saying that REM/MaxEnt can be viewed as an asymptotic instance of
the simple and self-evident Maximum Probability method.

Alternatively, \cite{Alberto} suggests to view REM/MaxEnt as a
separate method and hence to read the MaxProb/MaxEnt Thm as claiming
that REM/MaxEnt asymptotically coincides with MaxProb. If one adopts
this interesting and legitimate view then it is necessary to face
the fact that if $n$ is finite, the two methods in general differ.
{This} would open new questions. Among them also: MaxEnt/REM or
MaxProb? (i.e., $I$- or $\mu$-projection?) This is too delicate a
question to be answered by one sentence. Let us note, only, that
unless $n \rightarrow \infty$ entropy ignores multiplicity.

\subsection{$I$/$\mu$- or $\tau$-projection? MaxEnt/MaxProb or
MaxTent?}

The previous question (i.e., MaxEnt or MaxProb?) is a problem of
drawing an interpretational consequences from two variants of the
same probabilistic laws, and in this sense it can be viewed as an
'internal problem' of MaxEnt and MaxProb. From outside, from the
point of view of the Maximum R\'enyi-Tsallis Entropy method
(maxTent, \cite{Ts}, \cite{AJ}, \cite{Cohen})  MaxProb and MaxEnt
can be viewed as 'twins'.

maxTent is to the best of our knowledge intended by its proponents
for selection of probability distribution(s) under the setting of
BJIP with $\mathrm{\Pi}$ defined by $X$-frequency moment constraints
(cf. \cite{gg_nonlin}). It is not known whether such a feasible set
$\mathrm\Pi$ admits unique distribution with maximal value of
R\'enyi-Tsallis entropy (called $\tau$-projection at
\cite{gg_nonlin}) as it is also not known whether $I$-projection on
such a set is unique or not. The non-uniqueness makes  it difficult
to relay upon CWLLN when one wants to draw from an established
non-identity of $\tau$ and $I$-projection conclusion that maxTent
method violates CWLLN, cf. \cite{LS}. At \cite{gg_nonlin} this
difficulty has been avoided by considering an instance of the
$X$-frequency constraints where the feasible set reduced into a
convex set. Since $\tau$- and $I$-projection on the set were shown
to be different, CWLLN directly implies that maxTent in this case
selects asymptotically conditionally improbable distribution. The
Example below (taken from \cite{gg_nonlin}) illustrates the point.

\begin{ex} \cite{gg_nonlin} Let $\mathrm{\Pi} = \{p: \sum_{i=1}^3 p_i^2 (x_i -
b) = 0,\sum_{i=1}^3 p_i - 1 = 0\}$. Let  $\mathrm{X} = [-2\ \ 0\  \
1]$ and let $b = 0$. Then $\mathrm{\Pi} = \{p: p_3^2 = 2 p_1^2, \sum
p_i - 1 =0\}$ which effectively reduces to $\mathrm{\Pi} = \{p: p_2
= 1 - p_1(1 + \sqrt{2}), p_3 = \sqrt{2} p_1\}$. The source
 $q$ is assumed to be uniform $u$.

The feasible set $\mathrm{\Pi}$ is convex. Thus $I$-projection
$\hat{p}$ of $u$ on $\mathrm{\Pi}$ is unique, and can be found by
direct analytic maximization to be $\hat{p} = [0.2748 \ \ 0.3366 \ \
0.3886]$. Straightforward maximization of R\'enyi-Tsallis' entropy
leads to maxTent pmf $\hat{p}_T = [0.2735 \ \ 0.3398 \ \ 0.3867]$,
which is different than $\hat{p}$.

Convexity of the feasible set guarantees uniqueness of the
$I$-projection, and consequently allows to invoke CWLLN to claim
that any pmf from ${\mathrm\Pi}$ other than the $I$-projection has
asymptotically zero conditional probability that it will be
generated by $u$. 
\end{ex}

\smallskip

Obviously, $I$CET permits to show the fatal flow of maxTent in a
more direct and more general way.

\section{Further  Results}

Some further results related to asymptotic concentration of
conditional probability are contained in this Section.

\subsection{$r$-tuple $I$CET/CWLLN and MEM/GME Methods}

Maximum Entropy in the Mean method (MEM), or its discrete-case
relative, Generalized Maximum Entropy (GME) method, are interesting
extensions of the standard REM/MaxEnt method\footnote{For a tutorial
on MEM see \cite{Gzyl}. GME was introduced at \cite{GJM}, see also
\cite{MJM}.}. Though, usually a hierarchical structure of the
methods is highlighted, here a different feature of the method(s)
will appear to be important.

First, Golan-Judge-Miller ill-posed inverse problem (GJMIP) has to
be introduced. Its simple instance can be described as follows:  Let
there be two independent sources $q^1$, $q^2$ of sequences and hence
types. Let $\mathrm{X}$, $\mathrm{Y}$ be support of the first,
second source, respectively. Let a set $\mathrm{C}_n$ comprise  {\it
pairs} of the types $[\nu^{n,1}, \nu^{n,2}]$ which were drawn {\it
at the same time}. GJMIP amounts to selection of specific pair(s) of
types from $\mathrm{C}_n$ when the information $\{\mathrm{X},
\mathrm{Y}, n, q^1 \bot q^2, \mathrm{C}_n\}$ is supplied.

\begin{ex}  An example of GJMIP.  Let $\mathrm{X} = \mathrm{Y} = [1\ 2\ 3]$. Let
$q^1 = q^2 = [1/3 \ 1/3 \  1/3]$; $q^1 \bot q^2$; $(q^1 \mapsto
\nu^{n,1}) \wedge (q^2 \mapsto \nu^{n,2})$. Let $n = 100$,
$\mathrm{C}_n = \{[\nu^{n,1}, \nu^{n,2}]: \sum_{i=1}^3 \nu^{n,1}_i
x_i + \nu^{n,2}_i y_i = 4; \sum_{i=1}^3 \nu^{n,1}_i = 1;
\sum_{i=1}^3 \nu^{n,2}_i = 1 \}$.  Given this information, it is
necessary to select a pair (one or more) of types from
$\mathrm{C}_n$. 
\end{ex}

\smallskip

Since throughout the paper discrete and finite alphabet is assumed,
GME will be considered instead of MEM, in what follows. The
important feature of GME is that it selects jointly and
independently drawn pairs (or $r$-tuples) of types/pmfs. Thus, it is
suitable for application at the GJMIP context. An $r$-tuple
extension of CWLLN ($r$CWLLN) provides a probabilistic justification
to the GME, at the GJMIP context.

Given GJMIP information, GME selects from the feasible set of the
pairs of pmfs the one $[\hat{p}^1, \hat{p}^2]$ (or more) which
maximizes sum of the relative entropies with respect to $q^1$,
$q^2$; respectively.

\begin{rCWLLN}
Assume a GJMIP. Let $\mathrm{C}$ be convex, closed set. Let
$B([\hat{p}^1, \hat{p}^2], \epsilon)$ be an $\epsilon$-ball centered
at the pair
\begin{equation*}
[\hat{p}^1, \hat{p}^2] = \arg \sup_{[p^1, p^2] \in \mathrm{C}}
\sum_{i=1}^{m^1} p^1_i \log\frac{p^1_i}{q^1_i} + \sum_{j=1}^{m^2}
p^2_j \log\frac{p^2_j}{q^2_j}.
\end{equation*}
Let $n \rightarrow \infty$. Then
\begin{equation*}
\pi([\nu^{n,1}, \nu^{n,2}] \in B([\hat{p}^1, \hat{p}^2], \epsilon)
\, | \, [\nu^{n,1}, \nu^{n,2}] \in \mathrm{C}; (q^1 \mapsto
\nu^{n,1}) \wedge (q^2 \mapsto \nu^{n,2}); q^1 \bot q^2)  = 1
\end{equation*}
\end{rCWLLN}

Proof of $r$CWLLN can be constructed along the same lines as the
proof of CWLLN; the assumption that the pairs of sequnces/types are
drawn at the same time and from independent sources is crucial for
establishing the result. Similarly, $r$-generalization of $I$CET can
be formulated and proven; and obviously $\mu$-variants of the
results hold true.

$r$CWLLN permits to rise the same objections to application of
R\'enyi entropy based variant of GME in the GJMIP context, as those
that were risen to maxTent in the BJIP context.

Needless to say, $\mu$-variant of $r$CWLLN provides a probabilistic
justification  to MaxProb variant of GME.

\subsection{Bayesian Conditional Law of Large Numbers}

It is worth a brief mentioning,  that there is an inverse problem
which is in a sense antipodal to Boltzmann Jaynes Inverse Problem.
Let us call it the $\beta$-problem, after \cite{gj}.

One form of the $\beta$-problem can be formulated as follows: let
there be a set  $\mathrm Q$ of sources over which a prior
distribution $\pi(\cdot)$ is specified. Let  $\nu^n$ be an $n$-type
drawn from a source $r$, not necessarily in $\mathrm Q$. It is
necessary to select a source $q \in \mathrm Q$, given the
information-pentad $\{n, \mathrm X, \nu^n, \mathrm Q, \pi(\cdot)\}$.

Conditional Law of Large Number for Sources \cite{gL}, \cite{gPS} is
concerned about the asymptotic behavior of posterior probability
$\pi(q \in \mathrm B | (q \in \mathrm Q) \wedge \nu^n)$. It states
that, under certain conditions, the posterior probability
asymptotically piles up on the $L$-projection of $r$ on $\mathrm Q$.
Hence, the particular $\beta$-problem has to be solved by
$L$-divergence maximization method.

An application of Conditional Limit Theorem for Sources to criterion
choice problem associated with the Empirical Estimation \cite{MJM},
\cite{Owen} as well as further discussion, can be found at
\cite{gj}.

\section{Summary}

Conditional Equi-concentration of Types on $I$-projections -- an
extension of CWLLN to the case of non-unique $I$-projection -- was
presented. $I$CET states that the conditional concentration of types
happens on each of the proper $I$-projections in equal measure.
Also, Gibbs Conditioning Principle was enhanced to capture multiple
$I$-pro\-je\-ctions.  Extended  GCP says (when $t  = 1$) that
conditional probability of a letter is asymptotically given by the
equal-weight mixture of proper $I$-projection probabilities of the
letter. The conditional equi-concentration/equi-%
probability 'phenomenon' is in our view an interesting feature of
'randomness'. It might be of some interest also for Statistical
Mechanics as it resembles  phase coexistence of Thermodynamics (eg.
triple point of water, vapor and ice).

A general form of MaxProb/MaxEnt Thm, which states asymptotic
identity of $I$- and $\mu$-projections, was recalled. It permits to
formulate $\mu$-projection variants of the corresponding
$I$-projection laws: CWLLN/GCP/\-$I$CET/EGCP. In our view, the
$\mu$-variants allow for a deeper reading than their $I$-projection
counterparts -- since the $\mu$-laws express the asymptotic
conditional behavior of types in terms  of the most probable types.
For instance, $\mu$-projection variant of CWLLN says that types
conditionally concentrate on the asymptotically most probable one.
This is, in our view, more obvious statement than that made by
$I$-variant of CWLLN. MaxProb/MaxEnt Theorem is also instrumental
for establishing of $I$CET.

The main results -- Conditional Equi-concentration of Types (CET) in
both its $I$- and $\mu$-projection form as well as Extended Gibbs
Conditioning and Boltzmann Conditioning -- were supplemented also by
further considerations. They are summarized below.

Though $\mu$-projections and $I$-projections asymptotically
coincide, for a finite $n$  they are, in general, different. In
light of this fact, the asymptotic identity of $\mu$- and
$I$-projections can be viewed in two ways. Either as saying that 1)
$I$-projection of $q$ on $\mathrm\Pi$ is the asymptotic form of
$\mu$-projection of $q$ on $\mathrm{\Pi}_n$ or  that  2)
$\mu$-projections on $\mathrm{\Pi}_n$ and $I$-projections on
$\mathrm{\Pi}_n$ asymptotically coincide.  Regardless of the
preferred view, the $\mu$-variants of the laws provide a
probabilistic justification of Maximum Probability method (MaxProb,
cf.~\cite{gg_what}) (at least in the area of Boltzmann-Jaynes
inverse problem). If the second view is adopted, then, when $n$ is
finite, it is necessary to face the challenge of selecting between
REM/MaxEnt method and MaxProb method.

The results have a relevance also for Maximum R\'enyi-Tsallis
Entropy method (maxTent), which is over the last years in vogue in
Statistical Physics. maxTent is to the best of our knowledge
proposed as a method for solving BJIP, albeit with the feasible set
$\Pi$ defined by non-linear moment constraints. Since, in general,
maxTent distributions ($\tau$-projections on $\mathrm\Pi$) are
different than $I$/$\mu$-projections on $\mathrm\Pi$, $I$CET implies
that the maxTent method selects asymptotically conditionally
improbable/impossible distributions.

A straightforward extension of CWLLN/CET for $r$-tuples of types was
also mentioned. It was noted that the extension provides a
justification to the Generalized Maximum Entropy method in the area
of Golan-Judge-Miller Inverse Problem.

Conditional Law of Large Numbers for Sources and its implications
for the $\beta$-problem were also mentioned, in passing.


\section{Appendix}

\subsection{MaxProb/MaxEnt}

\begin{mm}
 Let $\mathrm{X}$ be finite set. Let $\mathrm{M}_n$ be
set of all $\mu$-projections of $q$ on $\mathrm{\Pi}_n$. Let
$\mathrm{I}$ be  set of all $I$-projections of $q$ on
$\mathrm{\Pi}$. For $n \rightarrow \infty$, $\mathrm{M}_n =
\mathrm{I}$.
\end{mm}

\begin{proof} \cite{gg_asy}
 Necessary and sufficient
conditions for $\hat{\nu}^n$ to be a $\mu$-projection of $q$ on
$\mathrm{\Pi}_n$ are: {\it a)} $\pi(\hat{\nu}^n; q) \ge \pi({\nu}^n;
q)$, $\forall \nu^n \in \mathrm{\Pi}_n$; {\it b)} whenever
$\tilde{\nu}^n$ has the property {\it a)} then $\pi(\hat{\nu}^n; q)
\le \pi(\tilde{\nu}^n; q)$. Requirement {\it a)} can be equivalently
stated as:
\begin{equation}
\left(\prod \frac{n_i!}{\hat{n}_i!}\right)^{1/n} \ge \left(\prod
q_i^{n_i - \hat{n}_i}\right)^{1/n}
\end{equation}
and {\it b)} similarly. Standard inequality $(n/e)^n < n! <
n(n/e)^n$ (valid for $n > 6$) allows to bind the LHS of (1):

\begin{equation}
 \frac{\prod (\nu^n_i)^{\nu^n_i}}{n^{m/n} \prod
(\hat{\nu}_i^n)^ {\hat{\nu}^n_i} (\prod \hat{\nu}_i^n)^{1/n}} <
\text{LHS} < \frac{n^{m/n} \prod ({\nu}_i^n)^ {{\nu}^n_i} (\prod
{\nu}_i^n)^{1/n}}{\prod (\hat{\nu}^n_i)^{\hat{\nu}^n_i}}
\end{equation}

and similar bounds can be stated in the case of the requirement {\it
b)}\footnote{Note that if an $i$-th  component $\nu_i^n$ of a type
is zero then it can be effectively omitted from calculations of
$\pi(\nu^n; q)$. Thus, it is  assumed that product operations at
(1), (2) are performed on non-zero components only.}. Since $m$ is
by assumption finite, as $n \rightarrow \infty$ the lower and upper
bounds at (2) collapse into the ratio $\prod
(\nu_i^n)^{(\nu_i^n)}/\prod (\hat{\nu}_i^n)^{(\hat{\nu}_i^n)}$.
Consequently, the necessary and sufficient conditions {\it a)}, {\it
b)} for $\mu$-projection turn as $n \rightarrow \infty$ into
(expressed in an equivalent log-form): {\it i)} $ \sum (\nu_i^n \log
\nu_i^n - \hat{\nu}_i^n \log \hat{\nu}_i^n) \ge \sum (\nu_i^n -
\hat{\nu}_i^n) \log q_i$ for all $\nu^n \in \mathrm{\Pi}_n$; and
{\it ii)} whenever $\tilde{\nu}^n$ has the property {\it i)} then
$\sum \hat{\nu}_i^n \log\hat{\nu}_i^n - \tilde{\nu}_i^n
\log\tilde{\nu}_i^n \ge \sum (\hat{\nu}_i^n - \tilde{\nu}_i^n)\log
q_i$.

Necessary and sufficient conditions for $\hat p$ to be an
$I$-projection of $q$ on $\mathrm\Pi$ are the following: {\it I)}
$\sum (p_i \log p_i - \hat{p}_i\log\hat{p}_i) \ge \sum (p_i -
\hat{p}_i)\log q_i$ for all $p \in \mathrm\Pi$; and {\it II)}
whenever $\tilde{p}$ has the property {\it I)} then $\sum \hat{p}_i
\log\hat{p}_i - \tilde{p}_i \log\tilde{p}_i \ge \sum (\hat{p}_i -
\tilde{p}_i)\log q_i$.

Comparison of {\it i)}, {\it ii)} and {\it I)}, {\it II)} then
completes the proof.
\end{proof}

\subsection{$I$CET}

The conditional equi-concentration of types can be seen as a
consequence of Sanov's Theorem and MaxProb/MaxEnt Theorem. Indeed,
Sanov's Theorem implies that the probability $\pi(\nu^n \in C; q)$
decays to zero for any open set $C$ which excludes all of the
$I$-projections. The asymptotic identity of $I$- and
$\mu$-projections shows that for $n \rightarrow \infty$, the
$I$-projections have the same value of the probability $\pi(\nu^n;
q)$.

The following is a rough attempt to make the argument a bit more
formal. It relays upon MaxProb/MaxEnt Thm and the Lemma, which
states a standard inequality for ratio of probabilities:

\begin{lem}
  Let $\nu^n$, $\dot{\nu}^n$ be two types from $\mathrm{\Pi}_n$.
Then
\begin{equation*}
\frac{\pi(\nu^n; q)}{\pi(\dot{\nu}^n; q)} <
\left(\frac{n}{m}\right)^m \prod_{i=1}^m \frac{
(\frac{q_i}{\nu_i^n})^{n \nu_i^n} }{ (\frac{q_i}{\dot{\nu}_i^n})^{n
\dot{\nu}_i^n}  }
\end{equation*}
\end{lem}

\begin{proof}
$\pi(\nu^n; q) \le \prod_{i=1}^m (\frac{q_i}{\nu_i^n})^{n \nu_i^n}$.
Since for $n > 6$, $(n/e)^n < n!$ $ < n(n/e)^n$, it follows that
$\pi(\dot{\nu}^n; q) > \frac{1}{\dot{n}_1 \dots \dot{n}_m}
\prod_{i=1}^m (\frac{q_i}{\dot{\nu}_i^n})^{n \dot{\nu}_i^n}$.
$\dot{n}_1 \dots \dot{n}_m < (\frac{n}{m})^m$.
\end{proof}

\begin{icet}
 Let $\mathrm{X}$ be finite. Let there be $\mathrm k$ proper
$I$-projections $\hat{p}^1, \hat{p}^2, \dots, \hat{p}^\mathrm{k}$ of
$q$ on $\mathrm\Pi$. Let $\epsilon > 0$ be such that for $j = 1, 2,
\dots, \mathrm{k}$ $\hat{p}^j$ is the only proper $I$-projection of
$q$ on $\mathrm\Pi$ in the ball $B(\hat{p}^j, \epsilon)$. Let $n
\rightarrow \infty$. Then for $j = 1, 2, \dots, \mathrm{k}$,
\begin{equation*}
\pi(\nu^n  \in B(\epsilon, \hat{p}^j) | \nu^n \in \mathrm{\Pi}; q) =
1/\mathrm{k}.
\end{equation*}
\end{icet}

\begin{proof} 
Clearly,
\begin{equation}
\pi(\nu^n \in B(\epsilon, \hat{p}^j) | \nu^n \in \mathrm{\Pi}; q
\mapsto \nu^n) \le \frac{ \sum_{\nu^n \in B} \pi(\nu^n;
q)}{\sum_{\nu^n \in \mathrm\Pi} \pi(\nu^n; q)}
\end{equation}

$B_n(\epsilon, \hat{p}^j) \triangleq B(\epsilon, \hat{p}^j) \cap
\mathrm{\Pi}_n$.

Without loss of generality, let there be unique $I$-projection
$\hat{p}^n_{B}$ of $q$ on the ball $B_n(\hat{p}^j, \epsilon)$.
(Sequence of the $I$-projections on $\mathrm{\Pi}_n$ converges to a
proper $I$-projection of $q$ on $\mathrm\Pi$. To an $I$-projection
on $\mathrm\Pi$ which is not proper, no sequence of
$I_n$-projections converges.) Also, without loss of generality let
there be $\mathrm{k}$ $I$-projections $\hat{p}^j_{\mathrm{\Pi}_n}$,
$j = 1, 2, \dots, \mathrm{k}$ of $q$ on $\mathrm{\Pi}_n$.

Let $\mathrm{A} \triangleq B_n\backslash\{\hat{p}^n_B\}$,
$\mathrm{B} \triangleq
\mathrm{\Pi}_n\backslash\{\hat{p}^j_{\mathrm{\Pi}_n}\}, j \neq 1$,
$\mathrm{C} \triangleq \mathrm{\Pi}_n\backslash \mathrm{B}$.

Then the Right-Hand Side of (3) can be rewritten as:

\begin{equation}
\frac{\pi(\hat{p}^n_B)}{\pi(\hat{p}^1_{\mathrm{\Pi}_n})}
\frac{1 + \frac{%
\sum_{\nu^n \in \mathrm{A}} \pi(\nu^n)} {\pi(\hat{p}^n_B)}}%
{1 + \frac{\sum_{\nu^n \in \mathrm{B}}
\pi(\nu^n)}{\pi(\hat{p}^1_{\mathrm{\Pi}_n})}  + \frac{\sum_{\nu^n
\in \mathrm{C}} \pi(\nu^n)}{\pi(\hat{p}^1_{\mathrm{\Pi}_n})}}
\end{equation}

By MaxProb/MaxEnt Thm $I$-projections have for $n \rightarrow
\infty$ the same and supremal value of $\pi(\cdot)$. This implies
that $\pi(\hat{p}^n_B)/\pi(\hat{p}^1_{\mathrm{\Pi}_n})$ converges to
$1$ (the case of $0/0$ limit is excluded by the supremity of
$\pi(\cdot)$). The same argument implies that the first ratio in the
denominator converges to $\mathrm{k} - 1$. The Lemma implies that
the ratio in the nominator as well as the second ratio in the
denominator converge to zero.
\end{proof}

\subsection{Extended GCP}

\begin{egcp}
 Let $\mathrm{X}$ be  a finite set.
 Let $\mathrm\Pi$ be such that it admits $\mathrm k$ proper $I$-projections $\hat{p}^1,
\hat{p}^2, \dots, \hat{p}^\mathrm{k}$ of $q$ on $\mathrm\Pi$. Then
for a fixed $t$,
\begin{equation*}
\lim_{n \rightarrow \infty} \pi(X_1 = x_1, \dots, X_t = x_t | \nu^n
\in \mathrm\Pi; q \mapsto \nu^n) = 1/\mathrm{k}
\sum_{j=1}^{\mathrm{k}} \prod_{l=1}^t \hat{p}_{x_l}^j.
\end{equation*}
\end{egcp}

\begin{proof}
Clearly,
\begin{equation}
\pi(X_1 = x_1, \dots, X_t = x_t | \nu^n \in \mathrm\Pi; q \mapsto
\nu^n) = \frac{\sum_{\nu^n \in \mathrm\Pi} \pi(X_1 = x_1, \dots, X_t
= x_t, \nu^n)}{\sum_{\nu^n \in \mathrm\Pi} \pi(\nu^n; q)}
\end{equation}

Let, in addition to partitioning used in proof of $I$CET,
$\mathrm{D} \triangleq \cup_{j=1}^{\mathrm k}
\{\hat{p}^j_{\mathrm{\Pi}_n}\}$.

Then the RHS of (5) can be rewritten as:

\begin{equation}
\frac{\sum_{\nu^n \in \mathrm{D}} \pi(X_1 = x_1, \dots, X_t = x_t,
\nu^n) + \sum_{\nu^n \in \mathrm{\Pi}_n\backslash\mathrm{D}} \pi(X_1
= x_1, \dots, X_t = x_t, \nu^n)}{\pi(\hat{p}^1_{\mathrm{\Pi}_n})(1 +
\frac{\sum_{\nu^n \in \mathrm{B}}
\pi(\nu^n)}{\pi(\hat{p}^1_{\mathrm{\Pi}_n})}  + \frac{\sum_{\nu^n
\in \mathrm{C}} \pi(\nu^n)}{\pi(\hat{p}^1_{\mathrm{\Pi}_n})})}
\end{equation}

MaxProb/MaxEnt Thm implies that the first ratio in the denominator
converges to $\mathrm{k} - 1$. By the Lemma, the second ratio in the
denominator of (6) converges to zero as $n$ goes to infinity. The
second term in the nominator as well goes to zero as $n \rightarrow
\infty$ (to see this, express the joint probability $\pi(X_1 = x_1,
\dots, X_t = x_t, \nu^n)$ as $\pi(X_1 = x_1, \dots, X_t =
x_t|\nu^n)\pi(\nu^n)$ and employ the Lemma).

Then, MaxProb/MaxEnt Thm implies, that  for $n \rightarrow \infty$
the RHS of (6) becomes equal to $1/\mathrm{k}
\sum_{j=1}^{\mathrm{k}} \pi(X_1 = x_1, \dots, X_t = x_t|\hat{p}^j)$.
Finally, invoke Csisz\'ar's 'urn argument' (cf. \cite{CsiMT}) to
conclude that the asymptotic form of the RHS of (6) is $1/\mathrm{k}
\sum_{j=1}^{\mathrm{k}} \prod_{l=1}^t \hat{p}_{X_l}^j$.
\end{proof}

\subsection{Rational $I$-projections}

Types can  concentrate, in some sense, on rational $I$-projection
$\hat{p} \in \mathrm{Q}^m$ even though the $I$-projection is
isolated point of the set $\mathrm{\Pi}$. The following Example
illustrates the concentration.

\begin{ex}
 Consider 
$\mathrm\Pi = \{p, \dot{p}\}$, where $p = [n_{1}/n_0, \dots,
n_{m}/n_0]$ and $\dot{p} = [\dot{n}_1/n_0,$ $\dots, \dot{n}_m/n_0]$,
$n_0 \in \mathrm{N}$. For $n \neq k n_0$, $k \in \mathrm{N}$ the set
$\mathrm{\Pi}_n$ is empty; otherwise it contains $p$ and $\dot{p}$.
In this case, concentration of types on $\mu$-projection is a direct
consequence of the next two Lemmas. The $I$-variant of the
concentration then arises from MaxProb/MaxEnt Thm.

\begin{lemN}
Let $\nu^n$, $\dot{\nu}^n$ be two n-types. Let $\delta = \nu^n -
\dot{\nu}^n$. Let $K$ denote the non-negative elements of $n\delta$,
$L$ the absolute value of negative elements of $n\delta$. Let $c =
\prod_{-} \dot{n}_i^{L_i}/\prod_{+} \dot{n}_i^{K_i}$, where the
subscript $-$, $+$ indicates that the index $i$ goes through the
elements of K, L, respectively. Then
$\frac{\Gamma(k\nu^{kn})}{\Gamma(k\dot{\nu}^{kn})} < c^k$, for any
$k \in \mathrm{N}$.
\end{lemN}

\begin{proof}
$\frac{\Gamma(k\nu^{kn})}{\Gamma(k\dot{\nu}^{kn})} = \frac{\prod_{-}
(k(\dot{n}_i - L_i) + 1)\cdots k\dot{n}_i}{\prod_{+} (k\dot{n}_i +
1)\cdots k(\dot{n}_i + K_i)}$. So
$\frac{\Gamma(k\nu^{kn})}{\Gamma(k\dot{\nu}^{kn})} \leq
\frac{\prod_{-} \dot{n}^{k L_i}_i}{\prod_{+} \dot{n}^{k K_i}_i}$,
which is just $c^k$.
\end{proof}

\begin{lemN}
Let types $\nu^n$, $\dot{\nu}^n$ be such that $\pi(\dot{\nu}^n \,;\,
q) < \pi(\nu^n \,;\, q)$. Then $\frac{\pi({k\nu^{kn}} \,;\,
q)}{\pi(k\dot{\nu}^{kn} \,;\, q)} \rightarrow 0$ as $k \rightarrow
\infty$.
\end{lemN}

\begin{proof}
By the assumption, $\prod q_i^{n_i - \dot{n}_i} <
\frac{\Gamma(\tilde{\nu}^n)}{\Gamma(\nu^n)}$. The gamma-ratio is, by
the Lemma 1, smaller or equal to $c$, as defined at the Lemma. Thus,
$\prod q_i^{n_i - \dot{n}_i} = \gamma c$, where $\gamma \in [0,1)$,
$\gamma \in \mathrm{R}$. By Lemma 1, for any $k \in \mathrm{Z}$,
$\frac{\pi(k{\nu}^{kn} \,;\, q)}{\pi(k\dot{\nu}^{kn} \,;\, q)} \leq
(1/c)^k \prod q_i^{k(n_i - \tilde{n}_i)}$. The RHS of the
inequality, $\gamma^k$, goes for $k \rightarrow \infty$ to zero,
which completes the proof.
\end{proof}

In this case, if $\mathrm{\Pi}$ admits several rational
$I$/$\mu$-projections, then clearly, types equi-con\-cen\-trate on
them.
\end{ex}

\label{lastpage-97}

\end{document}